\documentclass[12pt,twoside,reqno]{amsart}
\usepackage{epsfig,amsmath,amssymb,fancyhdr,url,graphicx,latexsym,multicol,multirow,algorithmic,algorithm}
\calclayout

\begin{document}
\leftline{ \scriptsize \it  Journal of Prime Research in Mathematics
Vol. {\bf 10}(2015), 106-121.}

\vspace{1.3cm}

\title
{A New Lanczos-type Algorithm for Systems of Linear Equations}

\author{Muhammad Farooq$^1$ and Abdellah Salhi$^2$}
\thanks{ {\enskip $^{1}$Department of Mathematics, University of Peshawar, Khyber Pakhtunkhwa, 25120, Pakistan. Email: mfarooq@upesh.edu.pk\\
$^{2}$Department of Mathematical Sciences, University of Essex, Wivenhoe Park, Colchester, CO4 3SQ, UK. E-mail: as@essex.ac.uk}}
\begin{abstract}
Lanczos-type algorithms are efficient and easy to implement. Unfortunately they breakdown frequently and well before convergence has been achieved. These algorithms are typically based on recurrence relations which involve formal orthogonal polynomials of low degree. In this paper, we consider a recurrence relation that has not been studied before and which involves a relatively higher degree polynomial. Interestingly, it leads to an algorithm that shows superior stability when compared to existing Lanczos-type algorithms. This new algorithm is derived and described. It is then compared to the best known algorithms of this type, namely $A_5/B_{10}$, $A_8/B_{10}$, as well as Arnoldi's algorithm, on a set of standard test problems. Numerical results are included.
\vskip 0.4 true cm
 \noindent
 {\it AMS SUBJECT} : Primary 65F10.\\
\end{abstract}
\maketitle
\bibliographystyle{plain}

\pagestyle{myheadings} \markboth{\centerline {\scriptsize Muhammad Farooq and Abdellah Salhi}}
         {\centerline {\scriptsize  A New Lanczos-type Algorithm for Systems of Linear Equations}}
\section{Introduction and background} \label{sec:intro}
The Lanczos algorithm, \cite{50:Lanczos,52:Lanczos,04:Broyden}, is an iterative process that has been primarily designed to calculate the eigenvalues of a matrix. However, it has found a wide application in the area of Systems of Linear Equations (SLE's) where it now is a well established solver. Its attraction resides in its efficiency as it only involves vector-to-vector and matrix-to-vector products. Moreover, in exact arithmetic, it converges to the exact solution in at most $n$ steps,  where $n$ is the dimension of the problem, \cite{52:Lanczos}. While efficiency is its strong point, stability is not. Indeed, it is well known to breakdown as orthogonality of the so called Lanczos vectors, generated during the solution process, is lost. Efforts to avoid this breakdown led to a flurry of papers particularly from Brezinski and his team, \cite{95:Baheux,80:Brezinski,93:Brezinski,94:Brezinski,92:Brezinski,99:Brezinski,00:Brezinski,02:Brezinski}, and others,
\cite{98:Bjorck,92:Zaglia,00:Calvetti,10:Farooq,12:Farooq,11:Salhi,97:Greenbaum,99:Guennouni,06:Gérard,79:Parlett,85:Parlett,87:Saad,87:Van,94:Ye}.

Several Lanczos-type algorithms have been designed and among them, the famous conjugate gradient algorithm of Hestenes and Stiefel, \cite{52:Hest}, when the matrix is Hermitian and the bi-conjugate gradient algorithm of Fletcher, \cite{76:Fletcher}, and the algorithm of Arnoldi, \cite{51:Arnoldi,03:Saad}, in the general case.

Lanczos-type algorithms are commonly derived from recurrence relations typically using Formal Orthogonal Polynomials (FOP's) of low degree, \cite{52:Lanczos,93:Brezinski,83:Draux,39:Szego}. Recurrence relations using relatively higher degree FOP's have not been investigated. Here, we set out to design an algorithm that is based on such recurrence relations and FOP's, and study its properties and, in particular, its stability.


\subsection{The Lanczos Process}
Consider the SLE
\begin{equation}
\textit{A}\textbf{x}=\textbf{b},
\end{equation}
\noindent where $\textit{A}\in \textit{R}^{n\times n}$, $\textbf{b}\in \textit{R}^{n}$ and $\textbf{x}\in \textit{R}^{n}.$

\smallskip
Let $\textbf{x}_{0}$ and $\textbf{y}$ be two arbitrary vectors in $\textit{R}^{n}$ such that $\textbf{y}\neq0$ then Lanczos method \cite{52:Lanczos} consists in constructing a sequence of vectors $\textbf{x}_{k}\in \textit{R}^n$ defined as follows
\begin{equation}\textbf{x}_{k}-\textbf{x}_{0}\in \textit{K}_{k}(\textit{A}, \textbf{r}_{0})=span(\textbf{r}_{0}, \textit{A}\textbf{r}_{0},\dots,\textit{A}^{k-1}\textbf{r}_{0}),\end{equation}
\begin{equation}\textbf{r}_{k}=\textbf{b}-\textit{A}\textbf{x}_{k}\bot\textit{L}_{k}(\textit{A}^{T}, \textbf{y})=span(\textbf{y}, \textit{A}^{T}\textbf{y},\dots,\textit{A}^{T^{k-1}}\textbf{y}),\end{equation}
\noindent where $\textit{A}^{T}$ denotes the transpose of $\textit{A}$.

Equation $(2)$ gives,
\begin{equation}
\textbf{x}_{k}-\textbf{x}_{0}=-\alpha_{1}\textbf{r}_{0}-\alpha_{2}\textit{A}\textbf{r}_{0}-\dots-\alpha_{k}\textit{A}^{k-1}\textbf{r}_{0}.
\end{equation}
\noindent Multiplying both sides by $\textit{A}$ and adding and subtracting $\textbf{b}$ on the left hand side gives
\begin{equation}
\textbf{r}_{k}=\textbf{r}_{0}+\alpha_{1}\textbf{r}_{0}+\alpha_{2}\textit{A}\textbf{r}_{0}+\dots+\alpha_{k}\textit{A}^{k-1}\textbf{r}_{0}.
\end{equation}
\noindent If we set
\[P_{k}(x)=1+\alpha_{1}x+...+\alpha_{k}x,\]
\noindent then we can write from $(5)$
\begin{equation}
\textbf{r}_{k}=P_{k}(\textit{A})\textbf{r}_{0}.
\end{equation}

\noindent From $(3)$, the orthogonality condition gives

\medskip
$({\textit{A}^T}^{i}\textbf{y},\textbf{r}_{k})=(\textbf{y}, \textit{A}^i\textbf{r}_{k})=(\textbf{y}, \textit{A}^i\textbf{P}_{k}(\textit{A})\textbf{r}_{0})=0$, for $i=0,...,k-1.$

\medskip
\noindent Thus, the coefficients $\alpha_{1}$,...,$\alpha_{k}$ form a solution of SLE's,
\begin{equation}
\alpha_{1}(\textbf{y}, \textit{A}^{i+1}\textbf{r}_{0})+...+\alpha_{k}(\textbf{y}, \textit{A}^{i+k}\textbf{r}_{0})=-(\textbf{y}, \textit{A}^{i}\textbf{r}_{0}), \mbox{ for } i=0,\dots,k-1.\end{equation}
\noindent If the determinant of the above system is not zero then its solution exists and allows to obtain $\textbf{x}_{k}$ and $\textbf{r}_{k}$. Obviously, in practice, solving the above system directly for increasing values of $k$ is not feasible; $k$ is the order of the iterate in the solution process.  We shall see now how to solve this system for increasing values of $k$ recursively, that is, if polynomials $P_{k}$ can be computed recursively. Such computation is feasible, since polynomials $P_{k}$ form a family of FOP's which will briefly be explained below.
\subsection{Formal Orthogonal Polynomials}
Define a linear functional $c$ on the space of reel polynomials by
\begin{center}
$c(x^i)=c_{i}$ \mbox{ for } $i=0,1,\dots$
\end{center}

\noindent where \begin{center} $c_{i}=({\textit{A}^T}^{i}\textbf{y},\textbf{r}_{k})=(\textbf{y}, \textit{A}^i\textbf{r}_{k})$ \mbox{ for } $i=0,1,\dots$\end{center}

\noindent Write the orthogonality condition as,
\begin{equation}
c(x^iP_{k})=0 \mbox{ for } i=0,\dots,k-1.
\end{equation}

\noindent The above condition shows that $P_{k}$ is the polynomial of degree at most $k$ and is a FOP with respect to the functional $c$, \cite{80:Brezinski}. The normalization condition for these polynomials is $P_{k}(0)=1$; $P_{k}$ exists and is unique if the following Hankel determinant
\[\textit{H}^{(1)}_{k}=\left\vert\begin{array}{cccc}
c_{1} & c_{2} & \cdots & c_{k}\\
c_{2} & c_{3} & \cdots & c_{k+1}\\
\vdots & \vdots &  & \vdots\\
c_{k} & c_{k+1} & \cdots & c_{2k-1}\\
\end{array}\right\vert\]\noindent is not zero. In that case we can write $P_{k}(x)$ as follows.

\begin{equation}
P_{k}(x)=\frac{\left\vert\begin{array}{cccc}
1 & x & \cdots & x^k\\
c_{0} & c_{1} & \cdots & c_{k}\\
\vdots & \vdots &  & \vdots\\
c_{k-1} & c_{k} & \cdots & c_{2k-1}\\
\end{array}\right\vert}{\left\vert\begin{array}{cccc}
c_{1} & \cdots & c_{k}\\
\vdots &  & \vdots\\
c_{k} & \cdots & c_{2k-1}\\
\end{array}\right\vert},
\end{equation}
\noindent where the denominator of this polynomial is $\textit{H}^{(1)}_{k}$, the determinant of the system (7). We assume that $\forall$ $k$, $\textit{H}^{(1)}_{k}\neq0$ and therefore all the polynomials $P_{k}$ exist for all $k$. If for some $k$, $\textit{H}^{(1)}_{k}=0$, then $P_{k}$ does not exist and breakdown occurs in the algorithm, \cite{93:Brezinski,92:Brezinski,00:Brezinski,02:Brezinski,92:Zaglia}.

\medskip
A Lanczos-type method consists in computing $P_{k}$ recursively, then $\textbf{r}_{k}$ and finally $\textbf{x}_{k}$ such that $\textbf{r}_{k}=\textbf{b}-\textit{A}\textbf{x}_{k}$, without inverting matrix $A$. In exact arithmetic, this gives the solution of the system $(1)$
in at most $n$ steps, where $n$ is the dimension of the SLE, \cite{93:Brezinski,99:Brezinski}.

\subsection{Notation and organization}
The notation introduced by Baheux, in \cite{95:Baheux,94:Baheux}, for recurrence relations with three terms is adopted here.
It puts recurrence relations involving FOP's $P_k(x)$ (the polynomials of degree at most $k$ with regard to the linear functional $c$)
and/or FOP's $P^{(1)}_k(x)$ (the polynomials of degree at most $k$ with regard to linear functional $c^{(1)}$, \cite{91:Brezinski})
into two groups: $A_i$ and $B_j$. Although relations $A_i$, when they exist, rarely lead to Lanczos-type algorithms on their own (the exceptions being $A_4$, \cite{95:Baheux,94:Baheux}, and $A_{12}$, \cite{10:Farooq}, so far), relations $B_j$ never lead to such algorithms for obvious reasons. It is the combination of
recurrence relations $A_i$ and $B_j$, denoted as $A_i/B_j$, when both exist, that lead to Lanczos-type algorithms.
In the following we will refer to algorithms by the relation(s) that lead to them. Hence, there are, potentially, algorithms $A_i$
and algorithms $A_i/B_j$, for some $i=1,2,\dots $ and some $j=1,2,\dots $

In this paper, a new algorithm based on a recurrence relation that has not been studied before, is derived. It is then compared to
three other algorithms, one of which is the Arnoldi algorithm.

\smallskip
The rest of the paper is organized as follows. In the next section, the Lanczos-type algorithm $A_8/B_{10}$, \cite{95:Baheux}, and the estimation of the coefficients of the recurrence relations $A_8$ and $B_{10}$ used to derive it are given. Section 3 is on the estimation of the coefficients of recurrence relation $A_{12}$, \cite{10:Farooq}, used to derive the new algorithm of the same name. Section 4 describes the test problems and reports numerical results. Section 5 is the conclusion and further work.

\section{Baheux algorithm $A_{8}/B_{10}$}
The choice of algorithm $A_8/B_{10}$, for comparison with our own is dictated by the fact that this is the most robust of the algorithms considered in \cite{95:Baheux,94:Baheux} on some of the problems considered here. So, outperforming this algorithm implies outperforming the rest of the algorithms considered therein.

For completeness, we recall the relevant relations between adjacent FOPs that lead to $A_8/B_{10}$ and their coefficients estimates. These are $A_8$ and $B_{10}$.
The details of algorithm $A_5/B_{10}$ are given in \cite{94:Baheux}.

\medskip
\subsection{Recurrence relation $A_{8}$}
Relation $A_{8}$ is
\begin{equation}\label{A8}
P_{k}(x)=(A_{k}x+B_{k})P^{(1)}_{k-1}(x)+(C_{k}x+D_{k})P_{k-1}(x),
\end{equation}
first investigated in \cite{95:Baheux,94:Baheux}. Its coefficients are estimated as

\begin{equation} B_{k}=0, \end{equation}
\begin{equation} C_{k}=0, \end{equation}
\begin{equation}D_{k}=1,\end{equation}

\noindent and

\begin{equation}\label{Ak}A_{k}=-\frac{c(x^{k-1}P_{k-1}(x))}{c(x^kP^{(1)}_{k-1}(x))}.\end{equation}

\noindent As we know
\begin{eqnarray}\label{C}\begin{cases}
c(x^kP_{k})=({\textit{A}^T}^k\textbf{y},P_{k}(\textit{A})\textbf{r}_{0})=(\textbf{y}_{k}, \textbf{r}_{k}),\\
c(x^kP^{(1)}_k)=({\textit{A}^T}^k\textbf{y},P^{(1)}_{k}(\textit{A})\textbf{r}_{0})=(\textbf{y}_{k}, \textbf{z}_{k}),\end{cases}
\end{eqnarray}

\noindent with $\textbf{y}_{k}=\textit{A}^T\textbf{y}_{k-1}$ and $\textbf{z}_{k}$ is defined in (\ref{zk}). Using (\ref{C}), equation (\ref{Ak}) becomes

\begin{equation}A_{k}=-\frac{(\textbf{y}_{k-1}, \textbf{r}_{k-1})}{(\textbf{y}_k, \textbf{z}_{k-1})}=-\frac{(\textbf{y}_{k-1}, \textbf{r}_{k-1})}{(\textbf{y}_{k-1}, \textit{A}\textbf{z}_{k-1})}.\end{equation}

\medskip
\subsection{Recurrence relation $B_{10}$}
This relation, first investigated in \cite{95:Baheux,94:Baheux}, is
\begin{equation}\label{B10}
P^{(1)}_{k}(x)=(A^{1}_{k}x+B^{1}_{k})P^{(1)}_{k-1}(x)+C^{1}_{k}P_{k}(x),
\end{equation}
\noindent Its coefficients are estimated as

\begin{equation}A^{1}_k=0, \end{equation}
\begin{equation}C^{1}_{k}a_k=1, \end{equation}

\noindent where $a_{k}$ is the coefficient of $x^k$ in $P_{k}(x)$ defined in (\ref{A8}) and

\begin{equation} B^{1}_k=-\frac{C^{1}_{k}c(\textbf{y}_{k}, \textbf{r}_{k})}{c(\textbf{y}_{k}, \textbf{z}_{k-1})}. \end{equation}


\noindent Equation (\ref{B10}) gives
\begin{equation}\label{zk}
\textbf{z}_{k}=B^{1}_{k}\textbf{z}_{k-1}+C^{1}_{k}\textbf{r}_{k}.
\end{equation}

\subsection{Algorithm $A_{8}/B_{10}$}

The pseudo-code of $A_{8}/B_{10}$, due to Baheux, \cite{95:Baheux,94:Baheux}, is as follows.

\begin{algorithm}
\caption{Algorithm $A_{8}/B_{10}$}
\begin{algorithmic}
\STATE Choose $\textbf{x}_{0}$ and $\textbf{y}$ such that $\textbf{y}\neq0$, and $\epsilon$ arbitrarily small and positive.\\
\STATE Set $\textbf{r}_0 = \textbf{b} - A\textbf{x}_0$,\\

\STATE $\textbf{z}_{0}=\textbf{r}_{0}$,\\

\STATE $\textbf{y}_{0}=\textbf{y}$,\\

\FOR {$k=0,1,2,\dots$,}

\STATE $A_{k+1}=-\frac{(\textbf{y}_{k}, \textbf{r}_{k})}{(\textbf{y}_{k}, \textit{A}\textbf{z}_{k})}$,\\

\STATE $\textbf{r}_{k+1}=\textbf{r}_{k}+A_{k+1}A\textbf{z}_{k}$,\\

\STATE $\textbf{x}_{k+1}=\textbf{x}_{k}-A_{k+1}\textbf{z}_{k}$.\\

\IF{$||\textbf{r}_{k+1}|| \geq \epsilon$,}

\STATE $\textbf{y}_{k+1}=A^{T}\textbf{y}_{k}$,\\

\STATE $C^1_{k+1}=\frac{1}{A_{k+1}}$,\\

\STATE $B^{1}_{k+1}=-\frac{C^{1}_{k+1}(\textbf{y}_{k+1}, \textbf{r}_{k+1})}{(\textbf{y}_{k}, A\textbf{z}_{k})}$,\\

\STATE $\textbf{z}_{k+1}=B^{1}_{k+1}\textbf{z}_{k}+C^{1}_{k+1}\textbf{r}_{k+1}$.
\ELSE \STATE Stop; solution found.\\
\ENDIF
\ENDFOR
\end{algorithmic}
\end{algorithm}

\pagebreak
\section{A new Lanczos-type algorithm}
In the following, a new recurrence relation which leads to a new variant of the Lanczos algorithm is considered.
\subsection{Recurrence relation $A_{12}$}
Consider the following recurrence relation, \cite{10:Farooq,12Farooq}
\begin{equation}\label{Pk}
P_{k}(x)=A_{k}[(x^{2}+B_{k}x+C_{k})P_{k-2}(x)+(D_{k}x^3+E_{k}x^2+F_{k}x+G_{k})P_{k-3}(x)],
\end{equation}
\noindent for $k \ge 3$, where $A_{k}$, $B_{k}$, $C_{k}$, $D_{k}$, $E_{k}$, $F_{k}$ and $G_{k}$ are constants to be determined using the normalization condition $P_{k}(0)=1$ and the orthogonality conditions $c(x^{i}P_{k})=0$, $\forall i=0,\dots,k-1$, $x^{i}$ being a monic polynomial of exact degree $i$. To find these coefficients, we proceed as follows. Since $\forall k$, $P_{k}(0)=1$, equation (\ref{Pk}) gives
\[1=A_{k}[C_{k}+G_{k}].\]
\noindent Multiplying both sides of (\ref{Pk}) by $x^{i}$ and then applying the linear functional $c$, we get

\begin{equation}\label{OC}
\begin{array}{l}
c(x^{i}P_{k})=A_{k}\{c(x^{i+2}P_{k-2})+B_{k}c(x^{i+1}P_{k-2})+C_{k}c(x^{i}P_{k-2})+D_{k}c(x^{i+3}P_{k-3})\\
+E_{k}c(x^{i+2}P_{k-3})+F_{k}c(x^{i+1}P_{k-3})+G_{k}c(x^{i}P_{k-3})\}.\\
\end{array}\end{equation}

\noindent Equation (\ref{OC}) is always true for $i=0,...,k-7$.

\noindent For $i=k-6$, we have
\[0=D_{k}c(x^{k-3}P_{k-3}).\]
\noindent Since $c(x^{k-3}P_{k-3})\neq0$, we have \[D_{k}=0.\]

\medskip\noindent For $i=k-5$, we get
\[0=E_{k}c(x^{k-3}P_{k-3}).\]
\noindent But $c(x^{k-3}P_{k-3})\neq0$; therefore
\[E_{k}=0.\]

\medskip\noindent For $i=k-4$, (\ref{OC}) gives
\begin{equation}\label{F}F_{k}=-\frac{c(x^{k-2}P_{k-2})}{c(x^{k-3}P_{k-3})}.\end{equation}

\medskip\noindent For $i=k-3$, $i=k-2$ and $i=k-1$ we get  the following equations respectively
\begin{equation}\label{a1}
B_{k}c(x^{k-2}P_{k-2})+G_{k}c(x^{k-3}P_{k-3})=-c(x^{k-1}P_{k-2})-F_{k}c(x^{k-2}P_{k-3}),
\end{equation}
\begin{equation}\label{a2}\begin{array}{l}
B_{k}c(x^{k-1}P_{k-2})+C_{k}c(x^{k-2}P_{k-2})+G_{k}c(x^{k-2}P_{k-3})\\=-c(x^{k}P_{k-2})-F_{k}c(x^{k-1}P_{k-3}),
\end{array}\end{equation}
\noindent and
\begin{equation}\label{a3}\begin{array}{l}
B_{k}c(x^{k}P_{k-2})+C_{k}c(x^{k-1}P_{k-2})+G_{k}c(x^{k-1}P_{k-3})\\=-c(x^{k+1}P_{k-2})-F_{k}c(x^{k}P_{k-3}).
\end{array}\end{equation}

Let $a_{11}$, $a_{12}$, $a_{13}$, $a_{21}$, $a_{22}$, $a_{23}$, $a_{31}$, $a_{32}$, and $a_{33}$ be the coefficients of $B_{k}$, $C_{k}$ and $G_{k}$ in equations (\ref{a1}), (\ref{a2}) and (\ref{a3}) respectively and let $b_{1}$, $b_{2}$ and $b_{3}$ be the corresponding right sides of these equations. If $\Delta_{k}$ represents the determinant of the coefficients matrix of the above mentioned system of equations then, we have

\medskip
\noindent $a_{11}=c(x^{k-2}P_{k-2})$, $a_{12}=0$, $a_{13}=c(x^{k-3}P_{k-3})$,

\medskip
\noindent $a_{21}=c(x^{k-1}P_{k-2})$, $a_{22}=c(x^{k-2}P_{k-2})$, $a_{23}=c(x^{k-2}P_{k-3})$,

\medskip
\noindent $a_{31}=c(x^{k}P_{k-2})$, $a_{32}=c(x^{k-1}P_{k-2})$, $a_{33}=c(x^{k-1}P_{k-3})$,

\medskip
\noindent $b_{1}=-c(x^{k-1}P_{k-2})-F_{k}c(x^{k-2}P_{k-3})=-a_{21}-F_{k}a_{23}$,

\medskip
\noindent $b_{2}=-c(x^{k}P_{k-2})-F_{k}c(x^{k-1}P_{k-3})=-a_{31}-F_{k}a_{33}$,

\medskip
\noindent $b_{3}=-c(x^{k+1}P_{k-2})-F_{k}c(x^{k}P_{k-3})=-s-F_{k}t$, where $s=c(x^{k+1}P_{k-2})$ and $t=c(x^{k}P_{k-3})$.

\medskip
\noindent Therefore, equations (\ref{a1}), (\ref{a2}) and (\ref{a3}) can be written as

\begin{equation}
a_{11}B_{k}+0C_{k}+a_{13}G_{k}=b_{1},
\end{equation}
\begin{equation}
a_{21}B_{k}+a_{22}C_{k}+a_{23}G_{k}=b_{2},
\end{equation}
\begin{equation}
a_{31}B_{k}+a_{32}C_{k}+a_{33}G_{k}=b_{3}.
\end{equation}
\noindent To solve for $B_{k}$, $C_{k}$ and $G_{k}$, Cramer's rule requires
\[\Delta_{k}=a_{11}(a_{22}a_{33}-a_{32}a_{23})+a_{13}(a_{21}a_{32}-a_{31}a_{22}).\]
\noindent If $\Delta_{k}\neq0$, then
\begin{equation}B_{k}=\frac{b_{1}(a_{22}a_{33}-a_{32}a_{23})+a_{13}(b_{2}a_{32}-b_{3}a_{22})}{\Delta_{k}},\end{equation}
\begin{equation}G_{k}=\frac{b_{1}-a_{11}B_{k}}{a_{13}},\end{equation}
\begin{equation}C_{k}=\frac{b_{2}-a_{21}B_{k}-a_{23}G_{k}}{a_{22}},\end{equation} \noindent and
\begin{equation}
1=A_{k}[C_{k}+G_{k}].
\end{equation}

\medskip
With all the necessary coefficients now determined, the expression of the polynomials $P_{k}(x)$ becomes
\begin{equation}\label{A12}P_{k}(x)=A_{k}\{(x^{2}+B_{k}x+C_{k})P_{k-2}(x)+(F_{k}x+G_{k})P_{k-3}(x)\}.\end{equation}

\medskip
\noindent Let us now use the relation $(\ref{A12})$ to compute $P_{k}(x)$, necessary for the computation of the residual $\textbf{r}_{k}=\textbf{b}-\textit{A}\textbf{x}_{k}=P_{k}(\textit{A})\textbf{r}_{0}$ and the corresponding vector $\textbf{x}_{k}$.

\medskip
Assume that $P_{k}$ has exact degree $k$ and the 3-term recurrence relationship (\ref{A12}) holds. To move to the Krylov space,
replace $x$ by $\textit{A}$ and multiply both side of $(\ref{A12})$ by $\textbf{r}_{0}$ to get,
\begin{equation}P_{k}(\textit{A})\textbf{r}_{0}=A_{k}[(\textit{A}^{2}+B_{k}\textit{A}+C_{k}\textit{I})P_{k-2}(\textit{A})\textbf{r}_{0}+(F_{k}\textit{A}+G_{k}\textit{I})P_{k-3}(\textit{A})\textbf{r}_{0}].\end{equation}

\noindent Using equation $(6)$, gives
\begin{equation}\label{rm}
\textbf{r}_{k}=A_{k}\{(\textit{A}^{2}+B_{k}\textit{A}+C_{k}\textit{I})\textbf{r}_{k-2}+(F_{k}\textit{A}+G_{k}\textit{I})\textbf{r}_{k-3}\}.
\end{equation}
\noindent And using $\textbf{r}_{k}=\textbf{b}-\textit{A}\textbf{x}_{k}$, gives
\begin{equation}\label{xm}
\textbf{x}_{k}=A_{k}\{C_{k}\textbf{x}_{k-2}+G_{k}\textbf{x}_{k-3}-(\textit{A}\textbf{r}_{k-2}+B_{k}\textbf{r}_{k-2}+F_{k}\textbf{r}_{k-3})\},
\end{equation}

\bigskip
\noindent with $F_{k}$ as in equation (\ref{F}). Using $(\ref{C})$, $F_{k}$ can be written as

\[F_{k}=-\frac{(\textbf{y}_{k-2}, \textbf{r}_{k-2})}{(\textbf{y}_{k-3}, \textbf{r}_{k-3})}.\]

\medskip
\noindent Condition $(\ref{C})$ can be used equally to rewrite the expressions of $a_{11}$, through $a_{33}$, $b_{1}$ to $b_{3}$ as follows.

\bigskip
\noindent $a_{11}=(\textbf{y}_{k-2}, \textbf{r}_{k-2})$, $a_{12}=0$, $a_{13}=(\textbf{y}_{k-3}, \textbf{r}_{k-3})$,

\bigskip
\noindent $a_{21}=(\textbf{y}_{k-1}, \textbf{r}_{k-2})$, $a_{22}=(\textbf{y}_{k-2}, \textbf{r}_{k-2})$, $a_{23}=(\textbf{y}_{k-2}, \textbf{r}_{k-3})$,

\bigskip
\noindent $a_{31}=(\textbf{y}_{k}, \textbf{r}_{k-2})$, $a_{32}=(\textbf{y}_{k-1}, \textbf{r}_{k-2})$, $a_{33}=(\textbf{y}_{k-1}, \textbf{r}_{k-3})$,

\bigskip
\noindent $b_{1}=-a_{21}-F_{k}a_{23}$, $b_{2}=-a_{31}-F_{k}a_{33}$, $b_{3}=-s-F_{k}t$,

\medskip
\noindent where $s=(\textbf{y}_{k+1},\textbf{r}_{k-2})$ and $t=(\textbf{y}_{k},\textbf{r}_{k-3})$.\\

\noindent These parameters allow the explicit computation of $B_k$, $G_k$, $C_k$, and $A_k$ as is given by equations 31, 32, 33 and 34 respectively. Equations (\ref{rm}) and (\ref{xm}) define the new Lanczos-type algorithm.

Now, since all previous formulae are only valid for $k\geq3$, it is necessary to find the expressions of the polynomials of degrees $1$ and $2$. From (9), we can write
\[P_{1}(x)=\frac{\left\vert\begin{array}{cc}
1 & x\\
c_{0} & c_{1}\\
\end{array}\right\vert}{c_{1}},\]
\[P_{1}(x)=1-\frac{c_{0}}{c_{1}}x,\]
$\textbf{r}_{1}=\textbf{r}_{0}-\frac{c_{0}}{c_{1}}\textit{A}\textbf{r}_{0}$, and $\textbf{x}_{1}=\textbf{x}_{0}+\frac{c_{0}}{c_{1}}\textbf{r}_{0}$,
\noindent where $c_{i}=(\textbf{y}, \textit{A}^i\textbf{r}_{0})$.\\\\
Using (9) again, we can write
\[P_{2}(x)=\frac{\left\vert\begin{array}{ccc}
1 & x & x^2\\
c_{0} & c_{1} & c_{2}\\
c_{1} & c_{2} & c_{3}\\
\end{array}\right\vert}{\left\vert\begin{array}{cccc}
c_{1} & c_{2}\\
c_{2} & c_{3}\\
\end{array}\right\vert},\]
\[P_{2}(x)=1-\frac{c_{0}c_{3}-c_{1}c_{2}}{c_{1}c_{3}-c_{2}^2}x+\frac{c_{0}c_{2}-c_{1}^2}{c_{1}c_{3}-c_{2}^2}x^2,\]
$\textbf{r}_{2}=\textbf{r}_{0}-\alpha\textit{A}\textbf{r}_{0}+\beta\textit{A}^2\textbf{r}_{0}$, and $\textbf{x}_{2}=\textbf{x}_{0}+\alpha\textbf{r}_{0}-\beta\textit{A}\textbf{r}_{0}$,
\noindent where $\alpha=\frac{c_{0}c_{3}-c_{1}c_{2}}{\delta}$, $\beta=\frac{c_{0}c_{2}-c_{1}^2}{\delta}$ and
$\delta=c_{1}c_{3}-c_{2}^2$.
\subsection{Algorithm $A_{12}$}
Putting together the various steps given in the above section, the new algorithm can be described as follows.
\begin{algorithm}[H]
\caption{Algorithm $A_{12}$}
\begin{algorithmic}
\STATE Choose $\textbf{x}_{0}$ and $\textbf{y}$ such that $\textbf{y}\neq0$, and choose $\epsilon$ arbitrarily small and positive.

\STATE Set $\textbf{r}_{0}=\textbf{b}-A\textbf{x}_{0}$, $\textbf{y}_{0}=\textbf{y}$, $\textbf{p}=A\textbf{r}_{0}$, $\textbf{p}_{1}=A\textbf{p}$,
$c_{0}=(\textbf{y}, \textbf{r}_{0})$, \STATE  $c_{1}=(\textbf{y}, \textbf{p})$, $c_{2}=(\textbf{y}, \textbf{p}_{1})$,
$c_{3}=(\textbf{y}, A\textbf{p}_{1})$, $\delta=c_{1}c_{3}-c_{2}^2$, \STATE
$\alpha=\frac{c_{0}c_{3}-c_{1}c_{2}}{\delta}$,
$\beta=\frac{c_{0}c_{2}-c_{1}^2}{\delta}$, \STATE
$\textbf{r}_{1}=\textbf{r}_{0}-\frac{c_{0}}{c_{1}}\textbf{p}$,
$\textbf{x}_{1}=\textbf{x}_{0}+\frac{c_0}{c_1}\textbf{r}_0$, \STATE
$\textbf{r}_{2}=\textbf{r}_{0}-\alpha p+\beta \textbf{p}_{1}$, $\textbf{x}_{2}=\textbf{x}_{0}+\alpha
\textbf{r}_{0}-\beta \textbf{p}$, \STATE $\textbf{y}_{1}=A^{T}\textbf{y}_{0}$, $\textbf{y}_{2}=A^{T}\textbf{y}_{1}$,
$\textbf{y}_{3}=A^{T}\textbf{y}_{2}$, $k=3$.

\WHILE {$||r_{k}||\ge \epsilon$}
\STATE $\textbf{y}_{k+1}=A^{T}\textbf{y}_{k}$,
$\textbf{q}_{1}=A\textbf{r}_{k-1}$, $\textbf{q}_{2}=A\textbf{q}_{1}$, $\textbf{q}_{3}=A\textbf{r}_{k-2}$, \STATE
$a_{11}=(\textbf{y}_{k-2}, \textbf{r}_{k-2})$, $a_{13}=(\textbf{y}_{k-3}, \textbf{r}_{k-3})$,
$a_{21}=(\textbf{y}_{k-1}, \textbf{r}_{k-2})$, $a_{22}=a_{11}$,\STATE
$a_{23}=(\textbf{y}_{k-2}, \textbf{r}_{k-3})$, $a_{31}=(\textbf{y}_{k}, \textbf{r}_{k-2})$,
$a_{32}=a_{21}$, $a_{33}=(\textbf{y}_{k-1}, \textbf{r}_{k-3})$, \STATE
$s=(\textbf{y}_{k+1}, \textbf{r}_{k-2})$, $t=(\textbf{y}_{k}, \textbf{r}_{k-3})$,
$F_{k}=-\frac{a_{11}}{a_{13}}$,\STATE
$b_{1}=-a_{21}-a_{23}F_{k}$, $b_{2}=-a_{31}-a_{33}F_{k}$,
$b_{3}=-s-tF_{k}$,\STATE
$\Delta_{k}=a_{11}(a_{22}a_{33}-a_{32}a_{23})+a_{13}(a_{21}a_{32}-a_{31}a_{22})$,
\STATE
$B_{k}=\frac{b_{1}(a_{22}a_{33}-a_{32}a_{23})+a_{13}(b_{2}a_{32}-b_{3}a_{22})}{\Delta_{k}}$,
\STATE $G_{k}=\frac{b_{1}-a_{11}B_{k}}{a_{13}}$, \STATE
$C_{k}=\frac{b_{2}-a_{21}B_{k}-a_{23}G_{k}}{a_{22}}$, \STATE
$A_{k}=\frac{1}{C_{k}+G_{k}}$, \STATE
$\textbf{r}_{k}=A_{k}\{\textbf{q}_{2}+B_{k}\textbf{q}_{1}+C_{k}\textbf{r}_{k-2}+F_{k}\textbf{q}_{3}+G_{k}\textbf{r}_{k-3}\}$,
\STATE
$\textbf{x}_{k}=A_{k}\{C_{k}\textbf{x}_{k-2}+G_{k}\textbf{x}_{k-3}-(\textbf{q}_{1}+B_{k}\textbf{r}_{k-2}+F_{k}\textbf{r}_{k-3})\}$,
\STATE $k=k+1$.
\ENDWHILE
\STATE Stop; solution found.\\
\end{algorithmic}
\end{algorithm}

\section{Numerical results}

$A_{12}$, \cite{10:Farooq}, the algorithm described in the above section, has been tested against algorithms $A_5/B_{10}$ and $A_8/B_{10}$, the best algorithms according to \cite{94:Baheux,95:Baheux}, as well as against the established Arnoldi algorithm, \cite{51:Arnoldi,03:Saad}.

\subsection{Test problems I}
The test problems considered here arise in the 5-point discretisation of the operator $\frac{-d^{2}}{dx^2}-\frac{d^{2}}{dy^2}+\gamma\frac{d}{dx}$ on a rectangular region. Comparative results on instances of the following problem ranging from dimension 10 to 100 for parameter $\delta$ taking values $0.0$ and $0.2$ respectively, are recorded in Table 1 and Table 2.

$A=\left(\begin{array}{ccccccc}
B & -I & \cdots & \cdots  & 0\\
-I & B & -I  &  & \vdots\\
\vdots & \ddots & \ddots & \ddots & \vdots\\
\vdots &  & -I & B & -I\\
0 & \cdots & \cdots & -I & B\\
\end{array}
\right),$ with
$B=\left(\begin{array}{ccccccc}
4 & \alpha & \cdots & \cdots & 0\\
\beta & 4 & \alpha &  & \vdots\\
\vdots & \ddots & \ddots & \ddots & \vdots\\
\vdots & & \beta & 4 & \alpha\\
0 & \cdots &  & \beta & 4\\
\end{array}\right).$\\

\noindent and $\alpha=-1+\delta$, and $\beta=-1-\delta$. The right-hand side is $\textbf{b}=\textit{A}\textbf{x}$,
where $\textbf{x}=(1, 1, \dots, 1)^{T}$, is the solution of the system. The dimension of $B$ is $10$.

\begin{table}[h]
\caption{Experimental results for problems when $\delta=0$}
\begin{center} 
\scalebox{0.75}{
\begin{tabular}{|c|c|c|c|c|c|c|c|c|c|c|}
\hline
 & \multicolumn{2}{c|}{Arnoldi} & \multicolumn{2}{c|}{$A_{5}/B_{10}$} &
\multicolumn{2}{c|}{$A_{8}/B_{10}$} & \multicolumn{2}{c|}{$A_{12}$}\\\cline{2-9}
$n$  & \multicolumn{1}{c}{$||r_{k}||$} & \multicolumn{1}{|c|}{t(sec)} &
\multicolumn{1}{c}{$||r_{k}||$} & \multicolumn{1}{|c|}{t(sec)} & \multicolumn{1}{c}{$||r_{k}||$} &
\multicolumn{1}{|c|}{t(sec)} & \multicolumn{1}{c}{$||r_{k}||$} & \multicolumn{1}{|c|}{t(sec)}\\\hline
$10$ & $1.2514E^{-10}$ & 0.001450 & $2.2940E^{-13}$ & 0.001892 & $1.7704E^{-13}$ & 0.002770 & $4.9623E^{-13}$ & 0.002819\\
$20$ & $1.7733E^{-11}$ & 0.002207 & $2.5256E^{-14}$ & 0.001842 & $1.7489E^{-13}$ & 0.002654 & $1.7536E^{-13}$ & 0.002904\\
$30$ & $1.2990E^{-14}$ & 0.003602 & $3.9026E^{-09}$ & 0.002220 & $4.9472E^{-09}$ & 0.003179 & $5.4705E^{-08}$ & 0.003370\\
$40$ & $3.5434E^{-11}$ & 0.006071 & $1.4770E^{-10}$ & 0.002416 & $8.4658E^{-10}$ & 0.003095 & $1.4776E^{-08}$ & 0.003526\\
$50$ & $6.1827E^{-08}$ & 0.008870 & $1.9959E^{-06}$ & 0.002962 & $1.3598E^{-06}$ & 0.003696 & $4.7994E^{-06}$ & 0.003980\\
$60$ & $2.9843E^{-14}$ & 0.012282 & $9.1910E^{-06}$ & 0.003001 & $3.7470E^{-06}$ & 0.003776 & $5.0010E^{-06}$ & 0.004354\\
$70$ & $4.2642E^{-13}$ & 0.017151 & $4.9035E^{-06}$ & 0.003622 & $4.2579E^{-06}$ & 0.004194 & $1.3781E^{-06}$ & 0.005658\\
$80$ & $5.0951E^{-08}$ & 0.021938 & $4.4311E^{-06}$ & 0.004498 & $7.7199E^{-06}$ & 0.005504 & $7.5581E^{-06}$ & 0.005271\\
$90$ & $9.6960E^{-13}$ & 0.029083 & $NaN$ &  & $8.5560E^{-06}$ & 0.007900 & $3.7301E^{-06}$ & 0.006541\\
$100$ & $1.1397E^{-13}$ & 0.036462 & $1.1889E^{-06}$ & 0.003849 & $3.1695E^{-06}$ & 0.004499 & $8.9530E^{-07}$ & 0.005084\\
\hline
\end{tabular}}
\end{center}
\label{turns}
\end{table}

\begin{table}[h]
\caption{Experimental results for problems when $\delta=0.2$}
\begin{center} 
\scalebox{0.75}{
\begin{tabular}{|c|c|c|c|c|c|c|c|c|c|c|}
\hline
 & \multicolumn{2}{c|}{Arnoldi} & \multicolumn{2}{c|}{$A_{5}/B_{10}$} &
\multicolumn{2}{c|}{$A_{8}/B_{10}$} & \multicolumn{2}{c|}{$A_{12}$}\\\cline{2-9}
$n$  & \multicolumn{1}{c}{$||r_{k}||$} & \multicolumn{1}{|c|}{t(sec)} &
\multicolumn{1}{c}{$||r_{k}||$} & \multicolumn{1}{|c|}{t(sec)} & \multicolumn{1}{c}{$||r_{k}||$} &
\multicolumn{1}{|c|}{t(sec)} & \multicolumn{1}{c}{$||r_{k}||$} & \multicolumn{1}{|c|}{t(sec)}\\\hline
$10$ & $2.3499E^{-15}$ & 0.001377 & $5.2347E^{-04}$ & 0.002339 & $5.2347E^{-04}$ & 0.002948 & $5.2347E^{-04}$ & 0.003231\\
$20$ & $5.6622E^{-11}$ & 0.002149 & $4.1778E^{-11}$ & 0.001842 & $5.8526E^{-11}$ & 0.003090 & $6.3915E^{-10}$ & 0.003372\\
$30$ & $6.8771E^{-15}$ & 0.003573 & $8.9881E^{-04}$ & 0.002220 & $8.9880E^{-04}$ & 0.003580 & $8.9880E^{-04}$ & 0.003847\\
$40$ & $1.8106E^{-10}$ & 0.006137 & $8.7583E^{-04}$ & 0.002830 & $9.3988E^{-04}$ & 0.003620 & $9.1261E^{-04}$ & 0.003977\\
$50$ & $3.5345E^{-08}$ & 0.008552 & $6.2669E^{-04}$ & 0.003360 & $5.7269E^{-04}$ & 0.004055 & $2.5040E^{-04}$ & 0.004964\\
$60$ & $2.8757E^{-13}$ & 0.012544 & $6.3670E^{-04}$ & 0.003877 & $8.4915E^{-04}$ & 0.004885 & $7.3489E^{-04}$ & 0.005345\\
$70$ & $4.2552E^{-13}$ & 0.017352 & $8.5670E^{-04}$ & 0.003902 & $7.0703E^{-04}$ & 0.006158 & $9.9086E^{-04}$ & 0.005052\\
$80$ & $1.7785E^{-04}$ & 0.021629 & $NaN$         &          & $NaN$         &          & $6.5602E^{-04}$ & 0.012131\\
$90$ & $1.4837E^{-04}$ & 0.029332 & $NaN$         &          & $7.5451E^{-04}$ & 0.011230 & $9.5294E^{-04}$ & 0.011842\\
$100$ & $5.8942E^{-13}$ & 0.037067 & $NaN$        &          & $NaN$         &          & $9.9710E^{-04}$ & 0.018899\\
\hline
\end{tabular}}
\end{center}
\label{turns}
\end{table}

\subsection{Test problems II}

The coefficient matrix here is taken as the Hilbert matrix, i.e. $\textit{A}=hilb(n)$, where $hilb(n)$ is a Matlab function, $n$ being the dimension of $\textit{A}$. The right-hand side $\textbf{b}$ and the solution $\textbf{x}$, are defined in the same way as in test problems I. The Hilbert matrix is notoriously ill-conditioned. Ill-conditioned systems of linear equations are notoriously difficult to solve to any useful accuracy, \cite{11Farooq,96:Kim,81:Rice}

\begin{table}[htp]
\caption{Experimental results when $A$ is a Hilbert matrix.}
\begin{center} 
\scalebox{0.75}{
\begin{tabular}{|c|c|c|c|c|c|c|c|c|c|c|}
\hline
   & \multicolumn{2}{c|}{Arnoldi} & \multicolumn{2}{c|}{$A_{5}/B_{10}$} &
\multicolumn{2}{c|}{$A_{8}/B_{10}$} & \multicolumn{2}{c|}{$A_{12}$}\\\cline{2-9}
$n$ &  \multicolumn{1}{c}{$||r_{k}||$} & \multicolumn{1}{|c|}{t(sec)} &
\multicolumn{1}{c}{$||r_{k}||$} & \multicolumn{1}{|c|}{t(sec)} & \multicolumn{1}{c}{$||r_{k}||$} &
\multicolumn{1}{|c|}{t(sec)} & \multicolumn{1}{c}{$||r_{k}||$} & \multicolumn{1}{|c|}{t(sec)}\\\hline
$10$ & $3.2101E^{-15}$ & 0.002364 & $4.4809E^{-06}$ & 0.001945 & $4.4809E^{-06}$ & 0.002834 & $4.4809E^{-06}$ & 0.002706\\
$20$ & $2.2288E^{-15}$ & 0.003285 & $9.6002E^{-05}$ & 0.001772 & $9.6002E^{-05}$ & 0.002839 & $9.6002E^{-05}$ & 0.002787\\
$30$ & $3.8953E^{-15}$ & 0.004148 & $1.3341E^{-05}$ & 0.001929 & $1.3376E^{-05}$ & 0.003050 & $1.3340E^{-05}$ & 0.002791\\
$40$ & $3.2251E^{-14}$ & 0.005521 & $3.8228E^{-05}$ & 0.001880 & $3.8276E^{-05}$ & 0.003327 & $3.8229E^{-05}$ & 0.002772\\
$50$ & $2.8673E^{-15}$ & 0.007068 & $7.9457E^{-05}$ & 0.001997 & $7.9463E^{-05}$ & 0.003475 & $7.9457E^{-05}$ & 0.003334\\
\hline
\end{tabular}}
\end{center}
\label{t3}
\end{table}

All algorithms have been implemented in Matlab version $7.8.0$ and run on a PC, under Microsoft Windows XP Professional Operating System, with 3.2GB RAM, and 2.40 GHz Intel(R) Core(TM) 2 CPU 6600. The problems are solved as dense problems, i.e. no sparsity has been exploited.
The results point to the Arnoldi algorithm being the most robust overall, but also the slowest overall. Algorithms $A_5/B_{10}$ and $A_8/B_{10}$ are the fastest overall, but the least robust overall; in fact they have failed to solve some problems in high dimension due to breakdown, of course, which is endemic in Lanczos-type algorithms. Algorithm $A_{12}$, like Arnoldi, solved all problems but faster and not as accurately. It is also more robust than $A_5/B_{10}$ and $A_8/B_{10}$ overall, but slower than both of them overall. Its lower speed compared to that of $A_5/B_{10}$ and $A_8/B_{10}$ is expected since the recurrence relation $A_{12}$ involves more coefficients
than both recurrence relations $A_5/B_{10}$ and $A_8/B_{10}$. Note that on the Hilbert-type problems, Table 3, the algorithms could not cope with dimensions higher than 50. Arnoldi is again the most stable overall and the slowest as the dimension grows. The other three algorithms have similar performances.

\section{Conclusion and further work}
The way Lanczos-type algorithms are derived using recurrence relations involving FOP's means that many such algorithms can be created, each based on a different set of relations. The choice of recurrence relations to use is dictated by the degree of FOP's to be involved; high degrees mean a large number of coefficients have to be calculated in the concerned Lanczos process.  This, consequently, dictates the computational complexity of the resulting Lanczos-type algorithm. However, it is well known, \cite{79:Khachyan,80:Lovasz}, that computational complexity does not always imply efficiency, or indeed robustness. Moreover, robustness and accuracy are often more important. It is therefore worthwhile to look beyond complexity issues sometimes, like we did here.

In this paper we have shown that, indeed, there are recurrence relations worth exploring since they lead to more robust algorithms. As a result, a new Lanczos-type algorithm, $A_{12}$ has been designed. The numerical performance of this algorithm is compared to that of two existing Lanczos-type algorithms, which were found to be the best among a number of Lanczos-type algorithms, \cite{95:Baheux,94:Baheux}, on the same set of problems as considered here. It is also compared to the well established Arnoldi algorithm. It is interesting to find that algorithm $A_{12}$ is overall more robust than $A_5/B_{10}$ and $A_8/B_{10}$ and faster than Arnoldi's. This makes it occupy, at least on the set of problems used here and elsewhere, a happy medium position. It is therefore the ideal candidate for time-limited applications which do not require high accuracy. 
\bibliography{lanczos1}
\end{document}